\documentclass[12pt,pdftex]{amsart}
\usepackage[pdftex]{graphicx}
\usepackage{amssymb}

\newtheorem{theo}{Theorem}[section]

\newcounter{deficislo}[section]
\newcommand{\prop}[1]{\refstepcounter{deficislo}{\noindent \bf
Proposition \thesection.\thedeficislo.\  }{\it #1}}

\newcommand{\pz}{P(z,\bar z)}

\newcommand{\ab}[1]{\vert z\vert^{#1}}
\newcommand{\cdva}{{\mathbb C^2}}

\newcommand{\te}{\theta}

\newcommand{\fz}{{F(z,\bar z,u)}}

\newcommand{ \al}{\alpha}
\newcommand{\pa}{\partial}

\begin{document}

\quad
\vskip 3.3cm

\title{
The local equivalence problem  in CR geometry}
\author[M. Kol\'{a}\v{r}]{Martin Kol\'a\v r}
\address{
Department of Mathematics, Masaryk University
{\it E-mail}: {\tt mkolar@math.muni.cz}}

\thanks{Supported by a grant of the GA CR no. 201/05/2117}

\begin{abstract}
This article  is dedicated to the centenary of the local CR equivalence
problem, formulated
 by Henri Poincar\'{e} in 1907. The first part  gives  an account of
Poincar\'{e}'s heuristic counting arguments,
suggesting existence of infinitely many local CR invariants.
 Then we sketch the beautiful completion of
Poincar\'{e}'s approach
to the problem
in the work of Chern and
Moser on Levi nondegenerate hypersurfaces. The last part is an
overview of recent
progress in solving the problem on
Levi degenerate manifolds.
\end{abstract}


\maketitle

\section{Introduction}

There are two fundamental facts which link  analysis and geometry
in one complex   variable. The local, almost obvious one, states
that every real analytic arc can be straightened  by an invertible
holomorphic map. The global one is the Riemann mapping theorem -
 any simply connected  subdomain of the complex plane
is biholomorphically equivalent to  the open unit disc.

In March of 1907, Rendiconti del Circolo Matematico di Palermo
published an article ``Les fonctions analytique de deux variables
et la repr\'{e}sentation conforme".
In the first section, titled ``\'{E}nnonc\'{e} du probl\`{e}m"
Poincar\'{e} asks the same questions in two complex variables. He
formulates the following problem, which was to become one of the
cornerstones of CR geometry:
\\[0.1mm]

\textsf{``Soit alors dans l'espace $zz'$ une portion de surface
\`a 3 dimensions $s$ et sur cette surface un point $m$. Soit dans
l'espace des $ZZ'$ une portion de surface \` a 3 dimensions $S$ et
sur cette surface un point $M$. Est-il possible de d\'{e}terminer
les fonctions $Z$ at $Z'$ de telle fa\c{c}on qu'elle soint
r\'{e}guli\`{e}re dans in voisinage du point $m$, que le point
$ZZ'$ soit en $M$ quand le point $zz'$ est en $m$ et qu'il
d\'{e}crive $S$ quand le point $zz'$ d\'{e}crit $s$? C'est le {\it
probl\`{e}m local.}"}
\\[1mm]

Poincar\'{e} then defines the global and the mixed version of the
problem. We will not attempt to survey all the development in
complex analysis inspired by the three equivalence problems.
Instead, this paper confines itself  to the local problem, and
results directly related to its solution in dimension two.

Since the notation and terminology,
being a century old, make the original article  less accessible today,
 we give in Section 2 a review of
its  two heuristic counting arguments, reformulated in modern
language. After introducing the first CR invariant, the Levi form,
 we sketch a solution of the local problem for Levi
nondegenerate hypersurfaces, obtained by S.\ S.\ Chern and J.\ K.\ Moser.
 Their
construction of normal forms  completes Poincar\'{e}'s approach
 for this class of manifolds.

In Section 5   we  consider Levi degenerate hypersurfaces and the
second  CR invariant - the type of the point - introduced by
J.\ J.\ Kohn in \cite{K}. Then we review the important result on
convergence of formal equivalences for hypersurfaces of finite
type, due to M.\ S.\ Baouendi, P.\ Ebenfelt and L.\ P.\ Rothschild.
Normal forms for finite type hypersurfaces are described in
Section 6. In Section 7 we consider  applications  to a
classification of local symmetry groups and  the jet determination
problem. Open problems, mainly for points of  infinite type,
are also discussed.

There are several excellent  surveys on closely related topics
in the literature. In particular, we mention the articles
\cite{BER}, \cite{IK}, \cite{We}. 

\section{Two counting arguments}

In more familiar notation, the local Poincar\'{e} problem asks if
for two given pieces of hypersurface $M_1$ and $M_2$ in
$\mathbb{C}^2$, and points $p_1\in M_1$ and $p_2 \in M_2$ there
exists a biholomorphic map in a neighbourhood of $p_1$ which maps
 $p_1$ to $p_2$ and $M_1$ to $M_2$.

Poincar\'{e} gives two heuristic counting arguments which suggest
that
the  problem does not always have a solution. In the first
one, he essentially derives the tangential Cauchy-Riemann
equation.

To review his argument,
 consider holomorphic coordinates $(z,w)$, where $z = x + iy$,
$w = u + i v$, and a biholomorphic transformation
\begin{equation}
z^* = f(z,w)\,,\quad \quad  w^* = g(z,w)\,,
\label{oni}\end{equation}
 where $f=f_1 + i f_2$ and $g = g_1 + i g_2$. The
components of $f$ and $g$ satisfy the ordinary Cauchy-Riemann
equations:
$$
\frac{df_1}{dx} = \frac{df_2}{dy}\,,\quad \quad  \frac{df_1}{dy}=
-\frac{df_2}{dx}
$$
and three other  such systems, one for $f$
replaced by $g$ and two other for derivatives with
 respect to $u$ and $v$.
That gives eight equations.

 Now
let $M_1$ be given by
 \begin{align}
v &= \Phi(x,y,u) \label{ona}\\
\intertext{and $M_2$ by}
v^*&= \Phi^*(x^*, y^*,
u^*)\,.\label{ono}
\end{align}
 Assuming that the
point $q$ in a neighbourhood of $p$ stays on the hypersurface
$M_1$, we  express all functions on $M_1$ in terms of the
three variables $x$, $y$, $u$. Differentiation with respect to these
three variables will be denoted by $\partial$, while
differentiation with respect to all four variables, considered
independent, will be denoted by ordinary $d$.  We have
$$
\frac{\pa f_1}{\pa x} = \frac{df_1}{dx} +
\frac{df_1}{dv}\frac{d\Phi}{dx}
$$
and eleven analogous equations
obtained by replacing $x$ by $y$ and $u$, and $f_1$ by $f_2$,
$g_1$ and $g_2$.

Now consider all the twenty equations obtained above and
eliminate the sixteen ordinary $d$-derivatives $\dfrac{df_1}{dx},
\dfrac{df_1}{dy}, \dots$. This leaves us with a system of four
linear partial differential equations for the twelve
$\pa$-derivatives  $\dfrac{\pa f_1}{\pa x}, \dots$, which we call
system S.

On the other hand, if the image of $q$ is to stay on $M_2$, we
have also
$$ \frac{\pa
g_2}{\pa x} = \frac{d\Phi^*}{dx^*}\frac{\pa f_1}{\pa x} +
\frac{d\Phi^*}{dy^*}\frac{\pa f_2}{\pa x} +
\frac{d\Phi^*}{du^*}\frac{\pa g_1}{\pa x}\,,
$$
and replacing $x$ by
$y$ and $u$ gives two other equations. Substituting these
expressions into S, we arrive at a system of four
differential equations for three unknown functions, $f_1, f_2,
g_1$, and their partial derivatives with respect to $x,y$ and $u$.
Hence, in general, it will be impossible to find a solution.

The second counting argument
 considers a refined
 version of the local equivalence problem. For given $ p_1 \in M_1$,  $p_2 \in M_2$
 and  $n \in \mathbb N$, we ask if there is a local biholomorphic map
 taking $p_1$ to $p_2$,  such that the image
 of $M_1$, denoted $M^{\prime}_1$, has $n$-th order of contact
  with $M_2$ at $p_2$. Without any loss of
 generality, assume $p_1 = p_2 =0$.

Again, let $M_1$ be given by a defining equation of the
 form
(\ref{ona}) and $M_2$ by  (\ref{ono}).  Consider the Taylor expansion of $\Phi^*$
 up to
order $n$. It involves
$$N  = \frac{(n+1)(n+2)(n+3)}{6} -1 $$
 arbitrary real
coefficients.
Next, consider a   transformation of the form (\ref{oni}).
Specifying $f$ and $g$ up to order $n$ involves
$2\big(\binom{n+2}{2}-1\big)$ complex coefficients, which is
$$
N^{\prime} = 2(n+1)(n+2) - 4 = 2n^2 + 6n$$ real coefficients.
Now we write the equation for the image of $M_1$  in a parametric
form with parameters $x$, $y$, $u$,
\begin{equation} x^* + i y^* = f\big(x,y,u, \Phi(x,y,u)\big)\,,\quad \quad
  u^*+ i v^* = g \big(x,y,u, \Phi(x,y,u)\big)\,.
 \end{equation}
In order to check if $M^{\prime}_1$ and $M_2$ have contact of order
$n$, we substitute those values  into the defining equation $v^* =
\Phi(x^*, y^*, u^*)$  for $M_2$. The order of contact is obtained
if all the resulting coefficients up to order $n$ are zero.
Considering the coefficients of $M_1$ and $M_2$ as given, and
those of $f$ and $g$ as unknown, we have $N$ equations for $N'$
unknowns. It remains to calculate  that $N
> N' $ if and only if
$$\frac{(n+1)(n+2)(n+3)}{6} -1 >2n^2 + 6n$$
which gives
$$n^2 > 6n + 25\,.
$$
The first integer which satisfies this relation is  $n = 9$.
Hence for $n \geq 9$ we cannot,
in general, get contact of order $n$.

 Poincar\'e also gives a heuristic definition of
local invariants. They are defined relative to a chosen reference
hypersurface, by considering infinitesimal perturbations of  the
hypersurface and of the identity transformation. In a loose sense,
the definition contains the ideas of using model hypersurfaces and
a certain linearization of the biholomorphism, which are central
for the construction of Chern and Moser,  and for later results on
Levi degenerate hypersurfaces.

 It is an interesting fact  that although Poincar\'{e}
gives this definition, he does not get to the point
of actually calculating  the lowest order invariant - the Levi
form. It was done two years later by E.\ E.\ Levi.

\section{The first invariant}

In fact, the first CR invariant appears  already at order two.
To define it, we now consider   a smooth hypersurface $M
\subseteq {\mathbb C}^n$, $n \geq 2$,
 and  a point
 $p\in M$.
 Let $r \in
C^{\infty}$ be a local defining function, i.e., for  a
neighbourhood $U$ of $p$
$$M \cap U = \big\{z \in U \mid  r(z) = 0\big\}\,,
$$
and $\nabla r \neq 0$ in $ M \cap U$.

The  Levi form on $M$ at $p$ is the Hermitian form defined by
$$
L_p(\zeta) = \sum_{i,j=1}^{n}\frac{\partial^2 r}{\partial z \partial
\bar z}(p)\zeta_i \bar \zeta_j\,,
$$
for $\zeta = (\zeta_1, \dots
\zeta_n)$ in the complex tangent space to $M$ at $p$,
$$
T_p^{\mathbb C}M = \Big\{ \zeta \in \mathbb C^n : \sum_{i=1}^n
\frac{\partial r}{\partial z_i}(p)\zeta_i = 0\Big\}\,.
$$

It is easily  checked the that the signature of the Levi form  is
a biholomorphic invariant of (an oriented) hypersurface.

In $\mathbb C^2$,
the complex tangent space is one dimensional, so the
Levi form is a scalar. Hence the first invariant of a hypersurface
 can be  thought of as taking values in the
three  point set $\{-1, 0, 1\}$. In the sequel, we consider this two
dimensional case.

 Not surprisingly, the equivalence problem
was first considered in the nondegenerate case, when the Levi form
is nonzero. The first substantial progress was made by  B.\ Segre
in 1931.
In \cite{S} he  defined a set of invariants, which
he thought   to form a complete set. A year later, E.\ Cartan showed
that in fact the set was  not complete, and provided himself a complete solution to
the problem, as an application of his general method of moving
frames. His intrinsic approach is  different from the
extrinsic one, using a reformulation of the problem in terms of
differential forms. For a detailed exposition of Cartan's solution
we refer the reader to the book of H.\ Jacobowitz (\cite{J}).

 The direct approach was again
taken up  in the first part of the celebrated paper of Chern and
Moser \cite{CM}. This part,  originating in the work of the second
author, solves the local Poincar\'{e} problem for Levi nondegenerate
 hypersurfaces, in arbitrary
dimension, by a construction of normal coordinates. 

\section{Chern-Moser normal form}

By his observations, Poincar\'{e} originated the extrinsic
approach to the problem, which directly analyzes the action of the
group of local biholomorphisms on the defining equation of the
hypersurface.

In this section  we consider  the case when the Levi form is
nondegenerate, i.e.\ the first invariant is nonzero,  and sketch
the solution of the Poincar\'e problem obtained in \cite{CM}.

We will use  again
local holomorphic coordinates $(z,w)$,
   centered at $p$, such that the hyperplane $\{ v=0
\}$ is tangent to $M$ at $p$. The complex tangent at $p$ is 
given by $\{ w = 0 \}$.

$M$ is locally described as a graph
of a function $v = \Phi(x,y,u)$.
Assuming  that $M$ is real analytic,
  $\Phi$ is the sum of its  Taylor expansion starting with 2-nd order
terms, which we will express in terms of $z$, $\bar z$, $u$.

The first step in normalizing $\Phi$ treats the leading second
order terms.
 We have
$$
v = \mbox{Re}\; \alpha z^2 + A \vert z \vert^2 + o\big(\ab 2, u\big)\,.
$$
By a change of variable $w^* = w +  \beta z^2$ we may eliminate
the harmonic term, taking $\beta = i\alpha $.
By definition, $A$ is the value of the Levi form at $p$,
corresponding to the defining function $ r = \Phi - v$, so $A \neq
0$.
 By a suitable scaling in the
$z$-variable and a change of sign in  $w$, if necessary,  we make
$A = 1$. Then we can write (with stars omitted)
\begin{equation}
v = \vert z \vert^2 + \fz\,,\label{bl}
\end{equation} where $F$ is
real analytic, with Taylor expansion
\begin{equation}\fz = \sum_{i+j+m\ge 2} a_{ijm} z^i \bar z^j
u^m\,,\label{on}
\end{equation}
where $a_{ijm}=\overline{a_{jim}}$
and $a_{110}= a_{200} = 0$.
 In the next step we will consider only transformations
 which preserve this form and normalize the higher order part
 $\fz$.

 The model hypersurface is defined using the leading term, as
$$
S = \big\{(z,w) \in \cdva  \mid  v = \ab 2 \big\}\,.
$$
$S$  is an unbounded  version of the unit sphere in $\mathbb
C^2$.

It  has a five dimensional group of local automorphisms,
consisting of transformations of the form
\begin{equation}z^* =
\frac{ \delta e^{i\theta} (z+ aw)}{(1 - 2i \bar a z -(\mu  + i
|a|^2)w)}\,, \quad \quad
  w^* =\frac{ \delta^2 w}{(1 - 2i \bar a z -(\mu  + i
  |a|^2)w)}\,,\label{gr}
\end{equation}
where $a\in \mathbb C$,  $\delta \in \mathbb R^*$ and  $\mu,
\theta \in \mathbb R$. We will denote this group by $\mathcal{H}
$.


One of the  ideas in [CM]
 is to
consider power series expansions  along real analytic curves
transversal to the complex tangent space at $p$, rather than the
ordinary expansion.  It reflects the inhomogeneity of the real
tangent space and the special role played by the transverse
coordinate $u$. Hence we will consider partial Taylor expansion of
$F$ in $z$, $\bar z$. Denoting
\begin{align*}
F_{ij}(u)& = \sum_{m=0}^{\infty} a_{ijm} u^m\,,
\\
\intertext{we  have}
\fz &= \sum_{i,j =0}^{\infty}F_{ij}(u)z^i\bar z^j\,.
\end{align*}

We will subject the defining equation to a general  biholomorphic
transformation
\begin{equation}
 z^*=z+ f(z,w)\,,\quad \quad
w^*=w+g(z,w)\,,
 \label{ja} \end{equation}
where $f$ and $g$ are
represented by power series
$$
 f(z,w)=\sum_{i,j =0}^{\infty} f_{ij} z^i w^j\,,\quad \quad
   g(z,w)= \sum_{i,j=0}^{\infty} g_{ij}
 z^i w^j\,.
$$
 The only requirement on (\ref{ja})
 is that it preserves form  (\ref{bl}).
Along with the partial Taylor expansion of $F$, we will consider
the corresponding expansions of $f$ and $g$. Denote
\begin{align*}
f_k(w) &= \sum_{j=0}^{\infty} f_{kj}w^j\,,
\quad \quad
 g_k(w) = \sum_{j=0}^{\infty} g_{kj}w^j\,,
\\
\intertext{so that}
f(z,w) &= \sum_{k=0}^{\infty}f_k(w)z^k\,,
\quad \quad  g(z,w) = \sum_{k=0}^{\infty}g_k(w)z^k\,.
\end{align*}

Now we can formulate  the normalizing conditions on $F$.


\begin{theo}[{\cite{CM}}] There exists a biholomorphic
change of coordinates such that the defining equation in the new
coordinates satisfies
\begin{equation}
\begin{array}{rl}
F_{j,0} = 0\,, &  \quad \quad j=0,1,\dots\,, \quad \quad \quad \quad \quad \quad \quad \quad  \\
F_{1,j} = 0\,,& \quad \quad j=1,2,3,\dots\,, \\
F_{2,2} = 0\,,& \quad \quad  \\
F_{3,3} = 0\,,& \\
F_{3,2} = 0\,.&
\end{array}
\label{form}
\end{equation}
 This transformation is
determined uniquely, up to a natural action of the symmetry group
$\mathcal{H} $.
\end{theo}

In order to give a few (heuristic) remarks about  the  proof,
consider the change of variables formula, obtained by substituting
(\ref{ja}) into $ \vert z^* \vert^2 + F^* = v^*$, and restricting
variables to $M$,
\begin{align*}
\big| z + f\big(z,u+i(\ab 2  + F)\big) \big|^2 &+ F^* \big( z+
f(z,u+i(\ab 2 + F)),  \overline{ z + f(\dots)}\,,\\[3pt]
\mbox{Re}\;
g(z,u+i(\ab 2 + F))\big ) =&\ \ab 2 + F + \mbox{Im}\ g\big(z,u+i(\ab 2
+ F)\big)\,,
\end{align*}
where the argument of $F$ is $(z,\bar z, u)$. It is viewed as an
equality of two power series in $z, \bar z , u$.  By multiplying
out we can in principle  obtain relations between various coefficients of $F^*$
and $F,f,g$.  Separating the leading linear term leads to the
Chern-Moser operator,
$$L(f,g) = \mbox{Re} \left(2\bar z f(z,u+i\ab 2)
+ i g( z,u+i\ab 2)\right)\,.$$

 The first two conditions in (\ref{form}) are relatively easy to
satisfy. Vanishing of the harmonic terms $F_{j,0}$ determines all
parts of $g$, except for $\mbox{Re}\; g_0$. Note that while $F_{j0}= 0$
for $j\geq 1$ is a complex condition, $F_{00}=0$ is a real
condition, which determines only one part of $g_0$, namely $\mbox{Im}\;
g_0$.

The coefficients  $F_{1j}$ of $z\bar z^{j} $ for $j \geq 2$ are
essentially absorbed into the leading term $z\bar z$ by the
substitution $z^* =  z + f_j(w)z^j$. This determines $f_j$ for all
$j\geq 2$.

In fact, in a more geometric setting, one can prove that for any
real analytic curve transversal to $T_p^{\mathbb C}M$ there is a
biholomorphic transformation which attains the first two
conditions and in the same time maps the given curve into the
u-axis. This curve can be chosen in such a way that $F_{32} = 0$,
which determines $f_0$. There is exactly one such curve in any
direction transverse to the complex tangent space at $p$. This
non-uniqueness corresponds to the parameter $a$ in (\ref{gr}).

The remaining  three conditions, $F_{11} = F_{22}=F_{33} =0$ then
determine $f_1$ and $\mbox{Re}\; g_0$. Geometrically, $F_{33} = 0$
corresponds to a choice of a preferred parametrization of the
curve. There is a projective one parameter family of such
parametrizations, corresponding to the parameter $\mu$ in
(\ref{gr}). Similarly, $f_1=0$ corresponds to a choice of a
preferred section of $T_{q}^{\mathbb C}M$ along the u-axis, which
is mapped into the unit section by the normalization mapping.
There is a unique such section for every initial condition given
by each vector in $T_p^{\mathbb C}M$. This non-uniqueness corresponds to the
parameters  $\delta$ and $\theta$ in (\ref{gr}).

\section{ Points of finite type }

When the Levi form vanishes, i.e.
$$L_p=0\,,$$
 one would like to find
  the next nontrivial invariant.
The second invariant, type of the point, was defined
 in the pioneering work of J.\ J.\ Kohn ([K]).

It can be defined as the maximal order of contact between complex
curves and $M$ at $p$
 (originally, it was  defined  in terms of commutators of CR vector
 fields, see \cite{K}, \cite{D}).

Note that  $M$ is Levi nondegenerate at $p$ if and only if $p$ is a point
of finite type two.
From now on we  assume  that $p$ is a point of finite type $k$,
where $k> 2$.

Since the structure of Levi degenerate points near a point of
finite type  is already quite diverse,
it doesn't seem reasonable to expect   any
kind of uniform geometric theory for such hypersurfaces. On the
other hand, the possibility of constructing a formal normal form theory
is still very attractive, since for many  applications convergence
is not necessary.

The first attempt to construct normal forms for Levi degenerate
hypersurfaces is due to P. Wong, who considered a class of
hypersurfaces of finite type four, given by
$$v = \ab 4 + a \ab 2 \mbox{Re}\; z^2 + \ab 2 u^2 + \dots\,,
 $$
where dots denote terms of order higher than four. Here $0 \leq  a
< \frac 43$, in particular $p$ is an isolated weakly pseudoconvex
point. His construction uses in an essential way both the leading
fourth order term in $z,\bar z$ and the additional  term
$\ab 2 u^2$ which controls the Levi form along the $u$-axis.
 Further results on the equivalence
problem and normal form constructions were obtained in  \cite{S},
\cite{E}, \cite{Ju}, \cite{BB}, \cite{BE}.


As a  first step  in  normalizing the defining equation $v =
\Phi(z, \bar z, u)$, we consider again the low order harmonic
terms.
One can  show easily  that $p\in M$ is a point of finite type $k$
if and only if
 there exist (uniquely determined) complex numbers
$\al_2, \dots, \al_{k} $ such that after the change of variable
$$w^* = w + \sum^{k}_{i=2} \al_i z^i
\ $$ the defining equation  has form
\begin{equation}v^* = \pz +
o\big(\ab k, u^* \big)\,,\label{po}
\end{equation}
where $P$ is a nonzero real valued homogeneous
polynomial of degree $k$ \begin{equation} \pz = \sum_{j=1}^{k-1}
a_j z^j\bar z^{k-j}\,.
 \label{onii}
\end{equation}
  Here $a_j \in
{\mathbb C}$ and $a_j = \overline{a_{k-j}}$, since $P$ is real valued.
Dropping stars we rewrite (\ref{po}) as
\begin{equation}
v = P(z,\bar z) + \fz \label{fz}\end{equation}
 and  define the homogeneous model to $M$ at $p$,
$$M_H = \big\{(z,w) \in \cdva\mid   v = \pz \big\}\,.
$$
$P$ is uniquely defined up to a linear change of variables

$$ w^* = \delta w \,, \quad \quad  z^* = \beta z\,,
 $$
where $\delta \in \mathbb R^* $ and $\beta \in \mathbb C^*$. The
subgroup of all such  transformations which preserve $M_H$
will be again denoted by $\mathcal{H}$.
 It
is straightforward  to determine this group explicitly. In most
cases $\mathcal{H}$ is the full local automorphism group of $M_H$
 (see
Section 7).


\section{Local equivalence of finite type hypesurfaces}

We will find a solution to the local equivalence problem for
finite type hypersurfaces
 by a
generalization of Chern-Moser's construction. Convergence of the normalizing
 map will not be proved.  In
fact it seems  plausible that it need not converge.
For the application we will need the essential result on
convergence of formal equivalences between finite type
hypersurfaces, due to Baouendi, Ebenfelt and Rothschild.

\begin{theo}[{\cite{BER}}]
{Let $M_1$, $M_2$ be two real
analytic hypersurfaces  in $\mathbb C^2$ and $p_1\in M_1$, $p_2
\in M_2$ be points of finite type. Let $\phi$ be a formal
equivalence between $(M_1, p_1)$ and $(M_2, p_2)$. Then $\phi$ is
convergent.}
\end{theo}

Proving this result involves intricate analysis of Segre
varieties.
We
refer the reader to \cite{BER} for a detailed description of this
technique.

 In order to prove that the  coefficients of the power series in
 normal form provide a complete set of
invariants, we have to show that two hypersurfaces which are
assigned the same power series
 are indeed
biholomorphically equivalent. The composition of the first
normalization mapping
 with
the inverse of the second one gives a formal equivalence of the
two hypersurfaces. The result of \cite{BER} implies convergence of
this formal equivalence.

In view of the convergence result, it is enough to find a
normalization on the level of formal power series.  The defining
function $\Phi$  and the transformation (\ref{ja}) will be
interpreted in this sense, the action of (\ref{ja}) on $\Phi$
being given  by the transformation rule (\ref{cov}) below.

Starting with  the finite type hypersurface (\ref{fz}), we give a
sketch of  the normal form construction. The most important
information carried by the model is its essential type, denoted by
$l$. It can be defined as the lowest index in (\ref{onii})  for
which $a_l \neq 0$, hence  $1\leq l \leq \frac{k}2$.

The equivalence problem now splits  into three cases,
depending on the form of the model.
 The most symmetric, circular  case, corresponding to $2l = k$ and
$P = a_l\ab k$. The tubular case,  when $P$ is equivalent to
$(\mbox{Re}\; z)^k$, which corresponds to a tube domain. All other
hypersurfaces  can be treated together, as the generic case.

In order to formulate the normal form conditions in the generic
case, we need  a natural
 scalar product on the vector space of homogeneous polynomials of degree $k-1$
without a harmonic term. If $Q = \sum^{k-2}_{j=1} \alpha_j z^j
\bar z^{k-1-j}$ and $S = \sum^{k-2}_{j=1} \beta_j z^j \bar
z^{k-1-j}$, we denote
$$(Q,S) = \sum^{k-2}_{j=1} \alpha_j \bar \beta_j\,.
$$
This notation is applied also to polynomials which contain a
harmonic term, which is  ignored. We extend this notation also
to polynomials whose coefficients are functions of $u$. In
particular, for $S = P_z = \sum_{j=1}^{k-2} j a_j z^{j-1}\bar
z^{k-j}$ we denote
$$(F_{k-1}, P_z) = \sum_{j=1}^{k-2}F_{j,k-1-j}  (j+1)\bar a_{j+1}\,.
$$
In the generic case we get the following normal form conditions.

\begin{theo}
[{\cite{Ko1}}]\label{theo7.2}
 {There exists a formal change of
coordinates such that  the new  defining equation
satisfies
\begin{align*}
F_{j0} & = 0\,, \ \ \ \ \ j=1,2,\dots,  \\
F_{k-l+j,l} & = 0\,, \ \ \ \ \ j= 1,2,\dots, \\
F_{k-l,l} & = 0\,, \\
F_{2k-2l, 2l} & = 0\,, \\
(F_{k-1}, P_z) & = 0\,.   \end{align*}
 It is determined uniquely up
to the action of the symmetry group $\mathcal H$.}
\end{theo}

 Since now
$\mathcal{H}$ contains only linear transformations, its action on
normal forms is straightforward.

We give again a few (heuristic) remarks about the proof. The first
two conditions are satisfied in a  similar way as in
 the Levi nondegenerate case. The first condition determines
$g_k$ for $k \geq 1$. The difference here is that we don't impose
the real condition $F_{00} = 0$ ( all the information
about $M$ which is used and all the normal form conditions are
complex). The second condition is again satisfied by absorbing the
terms $F_{k-l-j, l}$ into the leading term $z^{k-l} \bar z^{l}$,
which determines $f_k$ for $k\geq 2$.

 The third and fourth
condition, $F_{k-l,l} = F_{2k-2l, 2l} = 0$,  determine $f_1$ and
$g_0$. The scalar product condition  determines $f_0$.
In general,  there
seems to be no geometric interpretation of these conditions.

The construction in the circular case is similar to the
nondegenerate case. The model is now
\begin{equation}S_k = \{(z,w) \in \Bbb C^2 \ \vert \   v = \ab k
\}\,. \end{equation}

The local automorphism group of $S_k$ is three dimensional,
consisting of transformations of the form

\begin{equation}f( z,w)  =  \frac{ \delta e^{i\theta} z}{(1 + \mu w)^{\frac1l}}\,,
\ \ \ \ g(( z,w)=
  \frac{ \delta^k w}{1 + \mu w}\,,\label{4.1}\end{equation}
                           with $\delta > 0,$ and $\te, \mu  \in \Bbb R$.

 One obtains the following
normal form conditions:

\begin{align*}
F_{j0}&= 0\,, \ \ \ \ \ j=0,1,\dots,  \\
F_{l,l+j}&= 0\,, \ \ \ \ \ j= 0,1,2,\dots, \\
F_{2l, 2l}& = 0\,,\\
F_{3l, 3l}& = 0\,,   \\
F_{2l, 2l-1}& = 0\,. \ \end{align*}  and the same conclusion as in
Theorem \ref{theo7.2}. For normal forms in the tubular
case see \cite{Ko1}.


A fundamental tool for proving Theorem \ref{theo7.2} is again  the change
of variables formula
\begin{equation}  \Phi^*(z+f, \bar z+\bar f,u+\mbox{Re}\ g)  = \Phi(z, \bar z ,u) +
\mbox{Im}\ g\big(z,u+i\Phi(z, \bar z , u)\big)\,, \label{cov} \end{equation} where
$f$ and $\mbox{Re}\ g $ are also evaluated at $(z,u+i
 \Phi(z, \bar z , u))$.

 It becomes manageable if we assign weights to the variables,
namely weight one to $z, \bar z$ and weight $k$ to $u$.
  In this formula we separate the leading linear term. Denoting weights by subscripts,
for terms of weight $\mu>k$  we get
\begin{equation}\begin{aligned} 
\Phi^*_{\mu}(z, \bar z, u) & +2 \mbox{Re}\ P_z(z,\bar z)  f_{\mu-k+1}(z,u + i\pz) \\
 =&\ \Phi_{\mu}(z, \bar z, u) + \mbox{Im}\ g_{\mu}(z,u+i\pz) + \dots
 \end{aligned}\end{equation}
where dots denote terms depending on $f_{\nu-k+1}, g_{\nu},
F_{\nu}, F^*_{\nu}$ for  $\nu < \mu$ , and $P_z = \frac{\partial
P}{\partial z}.$ From this we obtain the generalized Chern-Moser
operator
$$L(f,g) =  \mbox{Re} \{ig(z,u+i\pz) +
2 P_z f(z,u + i\pz)\}\,.$$
 Careful analysis of this operator is an
essential part in the proof of Theorem \ref{theo7.2}.


\section{Applications and open problems}

The  normal form construction  gives immediately substantial information
about local automorphism groups and finite jet determination.
We will denote by  $\mbox{Aut}(M,p)$
the group of local automorphisms of $M$ at $p$
(i.e.\ local biholomorphic transformations  preserving $M$ and $p$).
By Theorem 6.2, the dimension of  $\mbox{\rm Aut}(M,p)$
is less then or equal to  the dimension of $\mbox{\rm Aut}(M_H,p)$.
In particular, in the generic case it is less or equal to one,
and less or equal to three in the circular case.


This result was sharpened in \cite{Ko2}, by further
analysis of the circular case
 and a refinement of the normal
forms, which leads to a full
classification of local symmetry groups for finite type
hypersurfaces.

\vskip 2mm
 \prop{ For a given hypersurface exactly one of
the following possibilities occurs.
\begin{enumerate}  \item $\mbox{\rm Aut}(M,p)$ has real dimension three. This
happens  if and only if $M$ is equivalent to $S_k$.

 \item $\mbox{\rm Aut}(M,p)$ has real dimension one and is noncompact,
  isomorphic to
$\Bbb R^+ \oplus {\Bbb Z}_m$.  This happens  if and only if $M$ is a
model hypersurface with $l<\frac k2$.

\item $\mbox{\rm Aut}(M,p)$ has real dimension one and is compact,  isomorphic to
$S^1$. This happens  if and only if  the defining equation of $M$ in
normal coordinates has form
$$ v = G(\ab 2, u).$$
\item $\mbox{\rm Aut}(M,p)$ is finite, isomorphic to $ {\Bbb Z}_m$. This
happens in all remaining cases.
\end{enumerate}
} The last case includes the trivial symmetry group, when $m=1$.
It is also  possible to determine the integer $m$ in the second and fourth
cases in terms of the
defining function in normal coordinates.

It can be seen from  (14) that   local automorphisms of $S_k$ are determined by their
2-jets.
 For all other hypersurfaces, 1-jets
are sufficient (\cite{Ko2}).

 \vskip 2mm
\prop  {\it Let $M$ be a hypersurface which is not equivalent to
$S_k$. Then local automorphisms are determined by their 1-jets.}

This result proves a conjecture formulated recently by Dmitri
Zaitsev  in the finite type case. His conjecture states that
1-jets suffice  for all hypersurfaces in dimension two, except for
those which are biholomorphic to the sphere at a generic point
(which is the case of $S_k$).

The local equivalence problem still remains open for points of
infinite type.
If $p$ is of infinite type, the order of contact of $M$ with
complex curves at $p$ is unbounded. It follows from real
analyticity that $M$ actually has to contain a complex
 curve. In the terminology of CR geometry, $M$ is not CR minimal,
 since it contains a proper submanifold of the same CR dimension
 as $M$, namely the complex curve.
In terms of local coordinates, $p$ is of infinite type if in suitable
  coordinates the hypersurface is given by
$$ v = u^s P(z,\bar z) + o(\ab k, u^s)\,.$$
Here $P$ is again a polynomial of degree $k$ of the form
(\ref{po}).  The numbers $s$ and $k$ are invariants of $M$.

 It is not known if  formal equivalences of such hypersurfaces are necessarily
 convergent. One additional difficulty is the fact that
 the 2-jet determination property for local automorphisms
 does not hold. For every integer $k$ there is an infinite type
 hypersurface whose local automorphisms are not determined
 by their k-jets (see \cite{Kow}, \cite{Z}).
In accord with Zaitsev's conjecture, all known examples are obtained
by blowing up the sphere.

On the other hand, there are some promising positive recent  results for certain classes of
 infinite type hypersurfaces (see e.g.\ \cite{ELZ}).

\vspace{5mm}

\vskip 1.5cm
\begin{thebibliography}{x--99}

\bibitem[BER]{BER}  Baouendi, M.\ S.,
Ebenfelt, P., Rothschild, L.\ P., \textit{Convergence and finite
determination of formal CR mappings},  J.\ Amer.\ Math.\ Soc.\
\textbf{13} (2000), 697--723.

\bibitem[BER2]{BER2} Baouendi, M.\ S.,
Ebenfelt, P., Rothschild, L.\ P.,   
 \textit{Local geometric properties of
real submanifolds in complex space}, Bull.\ Amer.\ Math.\ Soc.\ (N.S.) \textbf{37} 3
  (2000),   309--336.

\bibitem[BB]{BB}  Barletta, E., Bedford, E.,  \textit{Existence of proper
mappings from domains in $\Bbb C^2$ },  Indiana Univ.\  Math.\ J.\
 \textbf{2} (1990), 315--338.



\bibitem[BE]{BE} Beloshapka, V.\ K., Ezhov, V.\ V., \textit{Normal forms and 
model hypersurfaces in $\mathbb C^2$}, preprint.

\bibitem[BDS]{BDS} Burns, D., Jr.; Diederich, K.; Shnider, S.
\textit{ Distinguished curves in pseudoconvex boundaries,} 
 Duke Math. J.  \textbf{44}  (1977), no. 2, 407--431. 




\bibitem[CM]{CM} Chern, S.\ S.\ and Moser, J., \textit{Real hypersurfaces in
complex manifolds},  Acta Math.\  \textbf{133} (1974),  219--271.


\bibitem[D]{D} D'Angelo, J.\ P., \textit{ Orders od contact, real hypersurfaces and
applications},
 Ann.\ of Math. (2) \textbf{ 115} (1982), 615--637.


 \bibitem[E]{E}  Ebenfelt, P.,   \textit{New invariant tensors in CR
structures and a normal form for real hypersurfaces at a generic
Levi degeneracy}, J.\ Differential Geom.\  \textbf{50} (1998), 207--247.


\bibitem[ELZ]{ELZ} 
 Ebenfelt, P., Lamel, B., Zaitsev, D.,  \textit{
Degenerate real hypersurfaces in $\mathbb{C}^2$ with few
        automorphisms},  arXiv:math.CV/0605540


\bibitem[EHZ]{EHZ} Ebenfelt, P., Huang, X., Zaitsev, D.,
  \textit{The equivalence problem and rigidity for hypersurfaces embedded
  into hyperquadrics}, Amer.\ J.\ Math.\ \textbf{127} (2005), 169--191.




\bibitem[IK]{IK}  Isaev, A.\ V.,  Krantz,  S.\ G.,
 \textit{Domains with non-compact automorphism group: a survey.}
 Adv.\ Math.\  \textbf{146} 1  (1999),   1--38.


 \bibitem[J] {J}  Jacobowitz, H.,
\textit{An introduction to CR structures}, Math.\ Surveys
 Monogr.\  {\bf 32}, AMS 1990.


\bibitem[Ju]{Ju} Juhlin, R.,  \textit{PhD-thesis},  UCSD

 \bibitem[K]{K} Kohn, J.\ J., {Boundary behaviour of
$\bar \partial$ on weakly pseudoconvex manifolds of dimension two},
 J.\ Differential Geom.\  \textbf{6} (1972),  523--542.

\bibitem[Ko1]{Ko1} Kol\'a\v r, M.,
 \textit{Normal forms for hypersurfaces of finite type in
 $ \mathbb C^2$}, Math.\ Res.\ Lett.\ \textbf{12} (2005),  523--542.

\bibitem[Ko2]{Ko2} Kol\'a\v r, M., \textit{Local symmetries of  of finite type hypersurfaces in $\mathbb C^2$}, 
Sci.\ China A
  \textbf{48} (2006),  1633--1641.

\bibitem[Kow]{Kow} Kowalski, R., \textit{A hypersurface in $\Bbb C\sp 2$ whose
stability group is not determined by 2-jets},
  Proc.\ Amer.\ Math.\ Soc.\ \textbf{ 130} 12 (2002),  . 3679--3686.
 (electronic)



\bibitem[Me]{Me} Merker, J., \textit{
 Convergence of formal invertible CR mappings between minimal
 holomorphically nondegenerate real analytic hypersurfaces,} 
 Int. J. Math. Sci.  \textbf{26}  (2001),  no. 5, 281--302. 

\bibitem[Mi]{Mi} Mir,  N. \textit{
 Formal biholomorphic maps of real analytic hypersurfaces, }  
 Math. Res. Lett.  7  (2000),  \textbf{ 2-3}, 343--359. 

\bibitem[Po]{Po} Poincar\'e, H., \textit{Les fonctions analytique de
deux variables et la repr\'esentation conforme}, Rend.\ Circ.\ Mat.\
Palermo \textbf{23} (1907),  185--220.

\bibitem[S]{S} Segre, B.,
  \textit{Intorno al problem di Poincar\'e della rappresentazione pseudo-conform},
 Rend.\ Accad.\ Lincei \textbf{13} (1931), 676--683.

\bibitem[St]{St} Stanton, N., \textit{A normal form for rigid hypersurfaces in
$\Bbb C^2$ }, Amer.\ J.\ Math.\  \textbf{113} (1991),  877--910.

\bibitem[V]{V} Vitushkin, A.\ G., \textit{Real analytic
hypersurfaces in complex manifolds}, Russ.\ Math.\ Surv.\ \textbf{40}
(1985), 1--35.


\bibitem[We]{We} Wells, R.\ O., Jr., 
 \textit{The Cauchy-Riemann equations and differential geometry},
 Bull.\ Amer.\ Math.\ Soc.\ (N.S.)  \textbf{6} 2 (1982),  187--199.

\bibitem[Wo] {Wo}  Wong, P., \textit{A construction of normal forms for weakly
pseudoconvex CR manifolds in $\Bbb C^2$}, Invent.\ Math.\
\textbf{69} (1982),  311--329.

\bibitem[Z]{Z} Zaitsev, D.,
 \textit{Unique determination of local CR-maps by their jets: A survey},
  Atti  Accad.\ Naz.\ Lincei Cl.\ Sci.\ Fis.\ Mat.\ Natur.\ Rend.\ Lincei (9) Suppl.\ 
{\bf 13} (2002), 295--305.









\end{thebibliography}
\end{document}